\magnification=\magstep1

\centerline {\bf An explicit bound for the first sign change of the Fourier coefficients}

\centerline {\bf of a Siegel cusp form}
\bigskip\bigskip
\centerline {ŧ\it YoungJu Choie, Sanoli Gun and Winfried Kohnen}
\bigskip\bigskip\bigskip
\noindent {\bf Abstract  }
We give an explicit upper bound
for the first sign change of the Fourier coefficients
of an arbitrary non-zero Siegel cusp form $F$ of even integral weight on the Siegel modular group of arbitrary genus   $ g\geq 2 $.

\bigskip\bigskip

\noindent{\it{ Keynote}} :   Siegel cusp form,   Jacobi expansion,
Fourier expansion, sign change


\bigskip

\noindent {\bf 1. Introduction}
\bigskip
Fourier coefficients of cusp forms in general are quite mysterious objects. In particular, when real this applies to the distribution of their signs. Over the past years various aspects of the latter problem have been studied by several authors.\medskip
For example, in [9] it was shown  that the Fourier coefficients when real of a non-zero elliptic cusp form
$f$ on a congruence subgroup of the full modular group
$\Gamma_1:=SL_2({\bf Z})$ have infinitely many sign changes.
The proof which is rather  straightforward
uses the analytic properties of the Hecke $L$-function
and the Rankin-Selberg zeta function of $f$.
The above result was generalized in [8]
to the case of a non-zero Siegel cusp form $F$ of even integral weight on the symplectic group
$\Gamma_g:=Sp_g({\bf Z})\subset GL_{2g}({\bf Z})$
of arbitrary genus $g\geq 2$. In fact, in [8]
it was proved that if the Fourier coefficients $a(T)$ of $F$
are real ($T>0$ a positive definite half-integral matrix of size $g$),
then there exist infinitely many $T>0$ (modulo the usual action of $GL_g({\bf Z}))$
such that $a(T)>0$, and similarly such that $a(T)<0$. The proof uses similar arguments as for $g=1$,
with e.g. the Hecke $L$-function of $f$ replaced by the Koecher-Maass series of $F$.
\medskip

A more subtle problem is to give explicit upper bounds for the first sign change
in terms of the weight and the level.
If $f$ is a Hecke eigenform on the Hecke congruence group $\Gamma_0(N)$
of level $N$ this was done e.g. in [6,7,13].
The corresponding problem for arbitrary forms (not necessarily Hecke eigenforms)
 $f$ of squarefree level $N$ was studied in [3]. In fact, 
 in [3] the question was reduced to the case of Hecke eigenforms 
 by writing $f$ as a linear combination of the latter and 
 then using Chebyshev's inequality in combination with uniform lower bounds  
 for the Petersson norms of those Hecke eigenforms. 
 The technical details were rather complicated.
\medskip

In this paper we will give an explicit upper bound
for the first sign change of the Fourier coefficients
of an arbitrary non-zero Siegel cusp form $F$ of even integral weight 
on $\Gamma_g\, (g\geq 2)$. To the best of our knowledge, 
we think that this is the first general result in this direction.
The main idea will be to look at the Fourier-Jacobi expansion of $F$ where the coefficients are
Jacobi forms on the generalized Jacobi group
% $\Gamma_1\triangleright({\bf Z}^{g-1}\times {\bf Z}^{g-1})$ .
$\Gamma_1 \propto({\bf Z}^{g-1}\times {\bf Z}^{g-1})$ .
Using Taylor expansions of these coefficients we will reduce the question
to the case of elliptic modular forms and then will apply the results of [3].
Though this strategy appears rather simple, the actual technical details are somewhat involved.
\medskip

On the way we will also obtain a bound for the first non-vanishing Taylor coefficient function
of a generalized Jacobi form $\phi(\tau,z)\,(\tau\in {\cal H}$= upper half-plane, $z\in {\bf C}^{g-1}$)
around $z=0$. This generalizes a basic result of [4] in the case of 
classical Jacobi forms (Thm. 1.2) and eventually may be of independent interest.
\medskip
Although we do not have any direct and immediate application of our main result, 
we do think that the main steps in the proof may highlight again the theory
 of Jacobi forms, as an important bridge between Siegel modular forms 
 and classical elliptic modular forms.  
 Note that in the same spirit Jacobi forms played an important role 
beforehand in the proof of the Saito-Kurokawa conjecture [4] and
 in estimating Fourier coefficients of Siegel cusp forms [2]. 
\bigskip

\noindent {\bf 2. Statement of main result}
\bigskip

We will always suppose that $g\geq 2$. For $k\in {\bf N}$ we let $S_k(\Gamma_g)$
be the space of Siegel cusp forms
of weight $k$ on $\Gamma_g$.
Recall that this is the space of complex valued holomorphic functions $F(Z)$
on the Siegel upper half-space
${\cal H}_g$ (consisting of symmetric complex matrices of size $g$ with
positive definite imaginary part) such that
$$F((AZ+B)(CZ+D)^{-1})= \det (CZ+D)^kF(Z)$$
for all $\pmatrix {A&B\cr C&D\cr} \in \Gamma_g$ and such that $F$
has a Fourier expansion
$$F(Z)=\sum_{T>0}a(T)\, e^{2\pi i tr(TZ)} \qquad (Z\in {\cal H}_g),$$
where $T$ runs over all positive definite half-integral matrices of size $g$.
\medskip

We note that
$$a(T[U])=(\det U)^ka(T)$$
for all $T>0$ and all $U\in GL_g({\bf Z})$, where we have
used the standard abbreviation
$$A[B]:=B^tAB$$
for complex matrices $A$ and $B$ of appropriate sizes. \medskip In particular,
if $k$ is even and the $a(T)$ are real, then $a(T[U])$ is of the same sign as $a(T)$
 for all unimodular $U$.
\bigskip

\noindent {\bf Theorem.} {\it Let $k$ be even and $F\in S_k(\Gamma_g)$, $F\neq 0$.
Suppose that the Fourier coefficients $a(T)\, (T>0)$ of $F$ are real.
Then there exist $T_1>0, \; T_2>0$ with
$$tr\; T_1, \; tr\; T_2 \ll (k\cdot c_g)^5\, \log^{26}(k\cdot c_g)$$
such that $a(T_1)>0,\; a(T_2)<0$. Here
 $$c_g:= g\cdot 2^{g-1}\cdot (4/3)^{ g(g-1)/2}$$
% $$c_g:=  2^{g-1}\cdot (1+ {{g}\over{4\pi}} \cdot ({{4}\over{3}})^{  {g^2-g+1}\over{2} })$$
and the constant involved in $\ll$ is absolute and effective.}
\bigskip
\noindent {\it Remarks.} i) There is an obvious reformulation of 
the Theorem for $F\in S_k(\Gamma_g)$ a Siegel cusp form with arbitrary
 complex Fourier coefficients, with $a(T_1),\,a(T_2)$ replaced 
 by $Re (a(T_1)),\, \break Re (a(T_2))$  (and $Im (a(T_1)),\, Im (a(T_2))$). 
 This follows from the well-known fact that for $F\in S_k(\Gamma_g)$ 
 the Fourier series with general coefficients $Re (a(T))$ (resp. $Im (a(T))$) 
 are again in $S_k(\Gamma_g)$.
\medskip
ii) Any improvement of the bound obtained in [3] for elliptic cusp forms 
valid for all even integral weights would lead to a corresponding improvement
 of the bound in the Theorem, as will be clear from the proof.
\medskip
iii) One may ask more generally for the distribution of signs of 
the coefficients $a(T)$ where $T>0$ runs over primitive matrices only, 
or (oppositely) of the "radial'' coefficients $a(nT)\, ({n\geq 1})$, 
with $T>0$ fixed. The latter coefficients in the special case $g=2$ 
are related to the corresponding  Hecke eigenvalues (in case $F$ is a Hecke eigenform),
 see [1]. For some results regarding sign changes of eigenvalues we refer to [10,11,14]. 
 The method used here, however does not seem to give any new insights into the  questions addressed above.
\bigskip

The proof of the Theorem will be given in the next section.
\bigskip
\noindent {\bf 3. Proof}
\bigskip
It is sufficient to show the existence of $T>0$ with $tr\; T$
in the given range such that $a(T)<0$, since we can always replace $F$ by $-F$.
\medskip

Let us write
$$Z=\pmatrix{\tau & z\cr z^t &\tau'\cr}
\quad (\tau\in {\cal H},\, \tau'\in {ŧ\cal H}_{g-1}, \, z\in {\bf C}^{g-1}).$$
\medskip
Then $F$ has a Fourier-Jacobi expansion
$$F(Z)= \sum_{M>0}\phi_M(\tau,z) e^{2\pi i tr(M\tau')}\leqno (1)$$
where $M$ runs over all positive definite half-integral matrices of size $g-1$
 and the functions $\phi_M$ are Jacobi cusp forms of weight $k$
 and index $M$ on the generalized Jacobi group
$$\Gamma_1\propto({\bf Z}^{g-1}\times {\bf Z}^{g-1}).$$
\medskip

Those  are complex-valued holomorphic functions $\phi(\tau,z)$
on ${\cal H}\times {\bf C}^{g-1}$ with the transformation laws
$$\phi({{a\tau +b}\over {c\tau +d}}, {z\over {c\tau +d}})=(c\tau +d)^k\,
e^{2\pi i  {c M[z^t]}\over{ (c\tau+d)} } \,\phi(\tau,z)
 \quad (\forall \pmatrix {a&b\cr c&d\cr}\in \Gamma_1) \leqno (2)$$
and
$$\phi (\tau, z+\lambda\tau + \mu)= e^{-2\pi i  (M[\lambda^t] \tau   +2\lambda  z^t)}\,
\phi(\tau,z) \quad (\forall \lambda,\, \mu\in {\bf Z}^{g-1}),\leqno (3)$$
and having a Fourier expansion
$$\phi(\tau,z)= \sum_{n\geq 1,r\in {\bf Z}^{g-1}, 4n > M^{-1}[r^t]}c(n,r)\, e^{2\pi i(n\tau+rz^t)}$$
(we use the notation $A>B$ for symmetric real matrices $A$ and $B$ to indicate that $A-B>0$)(see [17]).
\medskip

Note that the $(n,r)$-th coefficient of $\phi_M(\tau,z)$ in (1) is equal to
$$a\Bigl(\pmatrix{n & r/2 \cr r^t/2 & M\cr}\Bigr)$$
(the condition $\pmatrix{n&r/2\cr r^t/2 & M\cr}>0$ is equivalent to
$n\geq 1,\,4n>M^{-1}[r^t]$ as follows from the usual Jacobi decomposition of the latter matrix).
\medskip

Since $F\neq 0$ there exists
$$T_0=\pmatrix{n_0 & r_0/2\cr r_0^t/2 & M_0\cr} >0$$
such that
$$tr\;T_0\ll_g k, \; a(T_0)\neq 0$$
as is well-known. In fact, one can find such a $T_0$
whose trace satisfies the explicit bound
$$tr \; T_0\leq {k\over {4\pi}}\cdot {2\over {\sqrt 3}}
\cdot g \cdot (4/3)^{g(g-1)/2} .\leqno (4).$$
As was kindly communicated to the authors by C. Poor, 
this is an easy  consequence of   basic reduction theory [5,12] and results proved in [15,16].
\medskip

Since $a(T_0)\neq 0$ the function $\phi_{M_0}(\tau,z)$
is not identically zero. We define
$$ \Phi_{M_0}(\tau,z):=\prod_\epsilon \phi_{M_0}(\tau,\epsilon z)$$
where $\epsilon=(\epsilon_1, \dots,\epsilon_{g-1})$ runs over all elements
of $\lbrace -	1,1\rbrace^{g-1}$ and
$\epsilon z:=(\epsilon_1 z_1, \dots, \epsilon_{g-1} z_{g-1}).$
\medskip

We  will simply write $\Phi(\tau,z)$ instead of  $\Phi_{M_0}(\tau,z)$.
\medskip

Note that
$$\phi_{M_0}(\tau,\epsilon z) = \phi_{M_0[D_\epsilon]}(\tau,z)$$
where $D_\epsilon$ is the diagonal matrix with entries
on the diagonal in the given order $\epsilon_1, \dots, \epsilon_{g-1}$.
This follows immediately if one acts on $F$ with the   matrix
$$\pmatrix{1 &0 & 0 & 0\cr 0 &D_\epsilon & 0 & 0 \cr 0 & 0 & 1 & 0 \cr 0 & 0 & 0 & D_{\epsilon}\cr}.$$
Therefore $\phi_{M_0}(\tau,\epsilon z)$ is a non-zero Jacobi cusp form
of weight $k$ and index $M_0[D_\epsilon]$
(which of course can also be checked by direct inspection).
\medskip

Hence we conclude that $\Phi(\tau, z)$ is a non-zero Jacobi cusp form
of weight $2^{g-1}k$ and index
$${\cal M}_0:=\sum_{\epsilon}M_0[D_\epsilon]. \leqno (5)$$
By construction this function is {\it even}  w.r.t. each of the variables
$z_\nu, \, \nu\in \lbrace 1, \dots, g-1\rbrace$. We observe that
$$tr\; {\cal M}_0 = 2^{g-1}\cdot tr\; M_0.$$
\medskip
We now develop $\Phi(\tau,z)$ in a Taylor series  around $z=0$, i.e. write
$$\Phi(\tau,z)=\sum_{\nu_1,\dots,\nu_{g-1}\geq 0}\chi_{\nu_1,\dots,
\nu_{g-1}}(\tau)z_1^{\nu_1}\dots z_{g-1}^{\nu_{g-1}}.$$
\medskip

Let us  choose $\alpha_{g-1}, \alpha_{g-2}, \dots, \alpha_2, \alpha_1$ in a ``minimal''
way such that $$\chi_{\alpha_1, \alpha_2, \dots, \alpha_{g-2},\alpha_{g-1}}(\tau)$$
is not the zero function. Here ``minimal'' means that for all $\tau$ we have
$$\chi_{\nu_1,\dots, \nu_{g-1}}(\tau)=0 \quad (0 \leq \forall \nu_{g-1}<\alpha_{g-1};
 \forall \nu_{g-2}, \dots, \nu_1\geq 0),$$
$$\chi_{\nu_1,\dots, \nu_{g-2}, \alpha_{g-1}}(\tau)
=0 \quad(0 \leq \forall \nu_{g-2}<\alpha_{g-2};
\forall \nu_{g-3}, \dots, \nu_1\geq 0),
 \dots \dots,$$
  $$\chi_{\nu_1, \alpha_2, \dots, \alpha_{g-1}}(\tau)
=0 \quad (0 \leq \forall \nu_1<\alpha_1).$$
\medskip

In the following, we will denote the diagonal elements of $M_0$ by
$$m_{11}^0, m_{22}^0,\dots, m_{g-1,g-1}^0~.$$
\bigskip
\noindent {\bf Lemma.} {\it The function
$$\sum_{\nu_1\geq 0}\chi_{\nu_1, \alpha_2, \dots,
\alpha_{g-1}}(\tau)z_1^{\nu_1}\leqno (6)$$
is a (classical) non-zero Jacobi cusp form of
weight $2^{g-1}k+\alpha_2+\dots +\alpha_{g-1}$
and index $2^{g-1}m_{11}^0$.}
\bigskip
\noindent {\bf Proof.} Up to a non-zero universal factor the function in (6) is equal to
$$\Bigl(\partial_{z_2}^{\alpha_2}\dots \partial_{z_{g-1}}^{\alpha_{g-1}}\,
\Phi(\tau, z_1, \dots, z_{g-1})\Bigr)_{|(z_2,\dots,z_{g-1})=(0,\dots, 0)}.$$
We differentiate  equation (2) successively w.r.t. $z_2, \dots, z_{g-1}$ up to
the orders $\alpha_2, \dots, \alpha_{g-1}$, respectively and use Leibniz
 rule together with the ``minimality'' of $\alpha_{g-1},\dots, \alpha_2$.
 Then we see that indeed (6) behaves like a Jacobi form of weight
 $2^{g-1}k+\alpha_2+\dots+\alpha_{g-1}$ and index $2^{g-1}m_{11}^0$ w.r.t.
 to the action of $\Gamma_1$.
\medskip

Likewise in (3) we take $\lambda=(\lambda_1, 0, \dots,0)$
and $\mu=(\mu_1,0, \dots, 0)$
and differentiate successively to see the correct behavior
of (6) under the action of ${\bf Z}^2$.
\medskip

The  conditions $n\geq 1,\, 4n>{\cal M}_0^{-1}[r^t]$
in the Fourier expansion of $\Phi$ are equivalent to
$$\pmatrix{n & r/2\cr r^t/2 & {\cal M}_0\cr} >0.$$
 Hence taking into account (5) (which implies that
 the $(1,1)$-entry of ${\cal M}_0$ is $2^{g-1}m_{11}^0$), it follows that
$$\pmatrix{n &r_1/2 \cr r_1/2 & 2^{g-1}m_{11}^0\cr} >0.$$
Therefore $4n\cdot 2^{g-1}m_{11}^0>r_1^2$ and
so the Fourier expansion of (6) is indeed as required.
\medskip

Finally we note that the function in (6) is not identically zero
since $\chi_{\alpha_1, \alpha_2, \dots, \alpha_{g-1}}(\tau)$ is
not the zero function by hypothesis. This proves the Lemma.
\bigskip

We now observe that $\chi_{\alpha_1, \dots,
\alpha_{g-1}}(\tau)$ is a non-zero cusp form of weight
$$k_1:= 2^{g-1} k+\alpha_1+\alpha_2+ \dots +\alpha_{g-1}\leqno (7)$$
on $\Gamma_1$, by the ``minimality'' of $\alpha_1$ and a similar argument
as above. If $g=2,$ full details are given in [4, p. 31]. 

  By Thm. 1.2, p. 10 in [4] it follows that
$$\alpha_1\ll2^{g-1}m_{11}^0.$$
\medskip

We now proceed inductively regarding the other variables $z_2, z_3 \dots$.
Thus we successively choose $\beta_{g-1},
\beta_{g-2}, \dots, \beta_3, \beta_1, \beta_2\geq 0$
in a ``minimal'' way such that
$\chi_{\beta_1,\beta_2, \dots, \beta_{g-1}}(\tau)$
is not the zero function, then successively
$\gamma_{g-1}, \gamma_{g-2}, \dots, \gamma_4, \gamma_2, \gamma_1, \gamma_3\geq 0$
in a ``minimal'' way such that
$\chi_{\gamma_1,\gamma_2, \gamma_3 \dots, \gamma_{g-1}}(\tau)$
 is not the zero function and so on.
\medskip

In the same way as above, using ``minimality'' one then shows that
$$\sum_{\nu_2\geq 0}\chi_{\beta_1, \nu_2, \beta_3 \dots, \beta_{g-1}}(\tau)z_2^{\nu_2}$$
is a non-zero Jacobi cusp form of weight
$2^{g-1}k+\beta_1+\beta_3+ \dots +\beta_{g-1}$
and index $2^{g-1}m_{22}^0$, the function
$$\sum_{\nu_3\geq 0}\chi_{\gamma_1, \gamma_2, \nu_3, \gamma_4,  \dots, \gamma_{g-1}}(\tau)z_3^{\nu_3}$$
is a non-zero Jacobi cusp form of weight
$2^{g-1}k+\gamma_1+\gamma_2+\gamma_4+  \dots +\gamma_{g-1}$
and index $2^{g-1}m_{33}^0$, and so on.
\medskip

Likewise as above, the functions $\chi_{\beta_1, \beta_2, \dots, \beta_{g-1}}(\tau)$ resp.
 $\chi_{\gamma_1,\gamma_2, \dots, \gamma_{g-1}}(\tau)$ and so
 on are non-zero cusp forms on $\Gamma_1$ of weights $2^{g-1}k+\beta_1+\dots +\beta_{g-1}$
resp. $2^{g-1}k+\gamma_1+\dots +\gamma_{g-1}$ and so on,
and $\beta_2\ll 2^{g-1}m_{22}^0, \,\gamma_3\ll 2^{g-1}m_{33}^0, \dots$.
\medskip

Since both the $\alpha$'s and the $\beta$'s are ``minimal'',
we conclude immediately from these conditions that
$$\beta_{g-1} =\alpha_{g-1}, \, \beta_{g-2}=\alpha_{g-2}, \dots, \beta_3=\alpha_3,$$
in a successive way. It then follows that
$$\alpha_2\leq \beta_2,$$
for otherwise we had
$$0=\chi_{\beta_1, \beta_2, \alpha_3, \dots,\alpha_{g-1}}(\tau)
= \chi_{\beta_1, \beta_2, \beta_3, \dots, \beta_{g-1}}(\tau),$$
a contradiction.
\medskip

We therefore conclude that
$$\alpha_2\ll 2^{g-1}m_{22}^0.$$
Proceeding in the same way with the $\alpha$'s and the $\gamma$'s, we infer that
$$\alpha_3\ll 2^{g-1}m_{33}^0.$$
Working on inductively, we finally conclude that
$$\alpha_1\ll 2^{g-1}m_{11}^0, \dots, \alpha_{g-1}\ll 2^{g-1}m_{g-1,g-1}^0,$$
hence
$$\alpha_1+\dots +\alpha_{g-1}\ll 2^{g-1} tr\; M_0.\leqno (8)$$
\bigskip

 We note that the arguments used above more generally
show  the following:
\bigskip

\noindent {\bf Proposition.} {\it  Let  $\phi(\tau,z)\, (\tau\in {\cal H}, z\in {\bf C}^{g-1})$ be
a generalized Jacobi form of weight $k$ and index $M>0$ on
$\Gamma_1\propto ({\bf Z}^{g-1}\times{\bf Z}^{g-1})$
and  let $\chi_{\nu_1, \dots,\nu_{g-1}}(\tau)\,
(\nu_1, \dots, \nu_{g-1}\geq 0)$ be  its Taylor coefficients around $z=0$.
Then there exists $(\alpha_1, \dots, \alpha_{g-1})\in
 {\bf N}_0^{g-1}$ such that $\alpha_1+ \dots +\alpha_{g-1}\ll tr\; M$
and $\chi_{\alpha_1, \dots, \alpha_{g-1}}(\tau)$ is a non-zero cusp form.}
\bigskip

This result generalizes Thm. 1.1 in [4]
in the classical case $\Gamma_1\propto {\bf Z}^2$ and may be of independent interest.
\bigskip

With the definition (7) it now follows from (8) that
$$k_1\ll 2^{g-1}(k+tr\; M_0). \leqno (9)$$
\medskip
Let us write $a(n)\, (n\geq 1)$ for the Fourier coefficients of $\chi_{\alpha_1,\dots,\alpha_{g-1}}(\tau)$.
Then by [3] (in the case of level 1) there exists $\tilde n\geq 1$ such that
 $$\tilde n \ll k_1^5 \,\log^{26}k_1, \,a(\tilde n)<0.$$
By (9) it follows that
$$\tilde n \ll \bigl(2^{g-1}(k+tr\;M_0)\bigr)^5\,\log^{26}\bigl(2^{g-1}(k+tr\;M_0)\bigr).\leqno (10)$$
\medskip

Observe that $\chi_{\alpha_1,\dots,\alpha_{g-1}}(\tau)$ up to a non-zero universal scalar equals
$$\Bigl(\partial_{z_1}^{\alpha_1}\dots \partial_{z_{g-1}}^{\alpha_{g-1}}\,\Phi(\tau, z)\Bigr)_{|z=0}.$$
Therefore if the Fourier coefficients of $\Phi(\tau,z)$ are
denoted by $C(n,r)$, it follows that up to a non-zero scalar
$a(\tilde n)$ is equal to
$$\sum_{r\in {\bf Z}^{g-1}, 4\tilde n >{\cal M}_0^{-1}[r^t]}
C(\tilde n,r)r_1^{ŧ\alpha_1}\dots r_{g-1}^{\alpha_{g-1}}.$$
Since $\Phi(\tau,z)$ is an even function w.r.t. each of
the variables $z_1, \dots, z_{g-1}$,  the $\alpha_1,
\dots, \alpha_{g-1}$ are all even integers. Hence it follows
that there exists $\tilde r\in {\bf Z}^{g-1}$
such that $C(\tilde n, \tilde r)<0$. However, $C(\tilde n, \tilde r)$
is a finite sum of products of Fourier coefficients
$$a\Bigl( \pmatrix{n_\epsilon & *\cr * & M_0[D_\epsilon]\cr}\Bigr)$$
where $\epsilon$ runs over $\lbrace -1,	1\rbrace^{g-1}$ as before, and with
$$\sum_\epsilon n_\epsilon = \tilde n.$$
Hence at least one of the coefficients
$$a\Bigl( \pmatrix{n_\epsilon & *\cr * & M_0[D_\epsilon]\cr}\Bigr)$$
must be negative.
\medskip
Since $n_\epsilon\leq \tilde n$ and $tr\; M_0 \leq tr\; T_0$ we infer from (10)
$$tr\; \pmatrix{n_\epsilon & *\cr * & M_0[D_\epsilon]\cr}=n_\epsilon + tr\;M_0$$
$$\ll  \bigl(2^{g-1}(k+tr\;M_0)\bigr)^5\,\log^{26}\bigl(2^{g-1}(k+tr\;M_0)\bigr) + tr\; T_0$$
$$\ll \bigl(2^{g-1}(k+tr\;T_0)\bigr)^5\,\log^{26}\bigl(2^{g-1}(k+tr\;T_0)\bigr).$$
Inserting from (4) we then obtain our assertion.
\bigskip

\noindent {\bf {Acknowledgements:}}  This work was partially supported by
 NRF 2013053914, NRF-2011-0008928 and NRF-2013R1A2A2A01068676.

\bigskip

\noindent {\bf References}
\bigskip
\noindent [1] $\;\;$ A.N.  Andrianov: Euler products corresponding to Siegel modular forms

of genus 2, Russian Math.  Surveys 29 (1974), 45-116.
\medskip
\noindent [2] $\;\;$ S. B\"ocherer and W. Kohnen: 
Estimates for Fourier coefficients of Siegel cusp forms, 

Math. Ann. 297 (1993), 499-517.
\medskip

\noindent [3] $\;\;$ Y. Choie and W. Kohnen:
The first sign change of Fourier coefficients of cusp forms,

Amer. J. Math. 131 (2009), no. 2,
517-543.
\medskip
\noindent [4] $\;\;$ M. Eichler and D. Zagier:
The theory of Jacobi forms. Progress in Math. vol 55,

Birkh\"auser 1985.
\medskip
\noindent [5] $\;\;$ E. Freitag:
Siegelsche Modulformen. Grundl. d. Math. Wiss. 254, Springer 1983.
\medskip
\noindent [6] $\;\;$ A. Ghosh and P. Sarnak:
Real zeros of holomorphic Hecke cusp forms,

J. Eur. Math. Soc. 14 (2012), 465-487.
\medskip
\noindent [7] $\;\;$ H. Iwaniec, W. Kohnen and J. Sengupta:
The first negative Hecke eigenvalue,

Int. J. Number Theory 3 (2007), 355-363.
\medskip
\noindent [8] $\;\;$ S. Jesgarz:
Vorzeichenwechsel von Fourierkoeffizienten von Siegelschen Spitzenformen,

Diploma Thesis (unpublished), Univ. of Heidelberg 2008.

\medskip
\noindent [9] $\;\;$ M. Knopp, W. Kohnen and W. Pribitkin:
On the signs of Fourier coefficients of

$\,$ cusp forms, Ramanujan J. 7 (2003), no. 1-3, 269-277.
\medskip
\noindent [10] $\;\;$ W. Kohnen: Sign changes of Hecke eigenvalues of Siegel

$\,$ cusp forms of genus two, Proc. Amer. Math. Soc. 135 (2007), 997-999.
\medskip
\noindent [11] $\;\;$ W. Kohnen and J. Sengupta: The first negative Hecke eigenvalue

$\,$ of a Siegel cusp form of genus two, Acta Arithmetica 129.1 (2007), 53-62.
\medskip
\noindent [12] $\;\;$ J. Martinet:
Perfect lattices in Euclidean spaces, Grundl. d. Math. Wiss. 327,

$\,$ Springer 2003.
\medskip
\noindent [13]  $\;\;$ K. Matom\"aki:
On signs of Fourier coefficients of cusp forms, Math. Proc. Camb.

$\,$ Phil. Soc. 152 (2012), 2007-2022.
\medskip
\noindent [14]  $\;\;$  A. Pitale and R. Schmidt: Sign changes of Hecke 
eigenvalues of Siegel cusp forms

$\;$of degree two, Proc. Amer. Math. Soc. 136 (2008), 3831-3838.
\medskip
 
\noindent [15] $\;\;$ C. Poor and D.S. Yuen:
Linear dependence among Siegel modular forms,

$\,$ Math. Ann. 318 (2000), 205-234.
\medskip
\noindent [16] $\;\;$ C. Poor and D.S. Yuen:
Dimensions of cusp forms for $\Gamma_0(p)$ in degree two

$\,$  and small weights, Abh. Math. Sem. Univ. Hamburg 77 (2007), 59-80.
\medskip

\noindent [17] $\;\;$ C. Ziegler:
Jacobi forms of higher degree,
  Abh. Math. Sem. Univ. Hamburg

  $\,$ 59 (1989), 191-224.

\bigskip\bigskip
{\it YoungJu Choie, Department of Mathematics, Pohang University of Science and

Technology, and PMI (Pohang Mathematical Institute), Pohang 790-784, Korea

E-mail: yjc@postech.ac.kr
\bigskip
Sanoli Gun, The Institute of Mathematical Sciences, C.I.T. Campus, Taramani,

Chennai 600 113, India

E-mail: sanoli@imsc.res.in
\bigskip
Winfried Kohnen, Mathematisches Institut der Universit\"at, INF 288,

D-69120 Heidelberg, Germany

E-mail: winfried@mathi.uni-heidelberg.de

\bye